\documentclass[a4paper,reqno,11pt]{amsart}

\usepackage{amssymb}
\usepackage{amsmath}
\usepackage{amsthm}
\usepackage{mathrsfs}
\usepackage{url}
\usepackage{color}
\usepackage{mathpazo}

\usepackage{pstricks}

\newtheorem{thm}{Theorem}

\newtheorem{pro}[thm]{Proposition}

\theoremstyle{remark}
\newtheorem{rem}[thm]{Remark}

\newtheorem{que}{Question}

\newtheorem*{ack}{Acknowledgement}

\def\D{\mathscr{D}}
\def\L{\mathscr{L}}
\def\R{\mathscr{R}}

\def\T{\mathcal{T}}

\def\ig#1{\mathsf{IG}(#1)}
\def\rig#1{\mathsf{RIG}(#1)}
\def\pre#1#2{\langle #1 \; | \; #2 \rangle}

\DeclareMathOperator\im{im}

\newcommand{\typeZ}{\boldsymbol{Z}}
\newcommand{\typeG}{\boldsymbol{G}}
\newcommand{\typeGbar}{\boldsymbol{\overline{G}}}
\newcommand{\typeR}{\boldsymbol{R}}

\newcommand{\rrel}{\boldsymbol{r}}

\psset{unit=1mm}
\newcommand{\incircle}[1]{\pspicture[shift=-1](0,0)(4,4)  \pscircle[fillstyle=solid,fillcolor=lightgray](2,2){2} \put(1.1,0.9){#1} \endpspicture}

\newcommand{\insquare}[1]{\pspicture[shift=-1](0,0)(4,4)  \psframe[fillstyle=solid,fillcolor=lightgray](0,0)(4,4) \put(1.1,0.9){#1} \endpspicture}

\newcommand{\indiamond}[1]{\pspicture[shift=-1](0,0)(4,4)  \psdiamond[fillstyle=solid,fillcolor=lightgray](2,2)(2,2.4) \put(1.1,0.9){#1} \endpspicture}

\begin{document}

\title[Free idempotent generated semigroups over bands]%
{Every group is a maximal subgroup of the free \\ idempotent generated semigroup over a band}

\author[I. Dolinka and N. Ru\v skuc]{Igor Dolinka and Nik Ru\v skuc}

\address{Department of Mathematics and Informatics, University of Novi Sad, Trg Do\-si\-te\-ja Obra\-do\-vi\-\'ca 4,
21101 Novi Sad, Serbia}

\email{dockie@dmi.uns.ac.rs}

\address{School of Mathematics and Statistics, University of St Andrews, St Andrews KY16 9SS, Scotland, UK}

\email{nik@mcs.st-and.ac.uk}



\dedicatory{Dedicated to Stuart W.\ Margolis on the occasion of his 60th birthday}

\begin{abstract}
Given an arbitrary group $G$ we construct a semigroup of idempotents (band) $B_G$ with the property that the free
idempotent generated semigroup over $B_G$ has a maximal subgroup isomorphic to $G$. If $G$ is finitely presented then
$B_G$ is finite. This answers several questions from recent papers in the area.
\end{abstract}

\maketitle

\section{Introduction}

Let $S$ be a semigroup. The set $E=E(S)$ of all idempotents of $S$ carries a structure of a partial algebra, called the
\emph{biordered set} of $S$, by retaining the products of the so-called \emph{basic pairs}: these are pairs of
idempotents $\{e,f\}$ such that $\{ef,fe\}\cap \{e,f\}\neq \emptyset$. It should be noted that if $ef\in\{e,f\}$ then
$fe$ is also an idempotent, possibly different from $e$, $f$ and $ef$. Also, if $S$ is an idempotent semigroup (i.e.\ a
\emph{band}) then its biordered set is in general different from $S$ itself, since not every pair is necessarily basic.
The term `biordered set' comes from an alternative (but equivalent) approach, where one considers $E(S)$ as a
relational structure equipped with two partial pre-orders; here we shall not pursue this approach, directing instead to
\cite{E1,E2,E3,Hi,Na} for further background.

The class of idempotent generated semigroups is of prime importance in semigroup theory, with a host of natural
examples, such as the semigroups of singular (non-bijective) transformations of a finite set (Howie \cite{Ho1}) or
singular $n\times n$ matrices over a field (Erdos \cite{Er}). It is not difficult to show that the category of all
idempotent generated semigroups with a fixed biordered set $E$ has an initial object $\ig{E}$, called the \emph{free
idempotent generated semigroup} over $E$ (we shall also say `over $S$' when $E=E(S)$). This semigroup is defined by the
presentation
$$\ig{E}=\pre{E}{e\cdot f=ef\ (\{e,f\}\text{ is a basic pair})}.$$
Here $e\cdot f$ stands for a word of length 2 in the free semigroup $E^+$, while $ef$ is the element of $E$ to which
the product equals in $S$. Unsurprisingly, $\ig{E}$ plays a crucial rule in understanding the structure of semigroups
with a prescribed biordered set of idempotents.

For reasons that are intrinsic to basic structure theory of semigroups \cite{Hi,Ho2}, this in turn  depends upon the
knowledge of maximal subgroups of $\ig{E}$. It was conjectured for a long time that the maximal subgroups of $\ig{E}$
are always free; this conjecture was widely circulated back in the 1980s, and was explicitly recorded in \cite{McE}.
The conjecture was proved in a number of particular cases, see e.g.\ \cite{McE,NP,Pa}. In 2009, Brittenham, Margolis
and Meakin \cite{BMM1} disproved the conjecture by means of an explicit  72-element semigroup $S$ such that $\ig{E(S)}$
has a maximal subgroup isomorphic to $\mathbb{Z}\oplus\mathbb{Z}$, the free abelian group of rank 2. This was followed
by Gray and Ru\v skuc \cite{GR1} who proved that \emph{every} group arises as a maximal subgroup of $\ig{E(S)}$ for a
suitably chosen semigroup $S$; if the group in question is finitely presented then a finite $S$ will suffice. Further
ensuing work such as \cite{GR2,DG, GY} investigates maximal subgroups of $\ig{S}$ for some specific natural semigroups
$S$, and the first author \cite{Do} initiates the study of $\ig{B}$, where $B$ is a band.

The aim of the present note is to prove the result announced in the title:

\begin{thm}\label{th1}
Let $G$ be a group. Then there exists a band $B_G$ such that $\ig{B_G}$ has a maximal subgroup isomorphic to $G$.
Furthermore, if $G$ is finitely presented, then $B_G$ can be constructed to be finite.
\end{thm}

This single construction provides an alternative, simpler proof of \emph{all} the main results of \cite{GR1} (Theorems
1--4), resolves \cite[Problem 1]{GR1} which asks whether every finitely presented group is a maximal subgroup of
$\ig{S}$ for some finite \emph{regular} semigroup $S$, and solves \cite[Problem 2]{Do} which calls for a
characterisation of maximal semigroups of free idempotent generated semigroups over bands.

\section{Presentation for maximal subgroups}

A general presentation for maximal subgroups of $\ig{S}$ in terms of parameters that depend only on the structure of
$S$ has been exhibited in \cite[Theorem 5]{GR1}. Since we are interested here only in the case of bands, we utilise the
particular form of this theorem, deduced in \cite[Corollary 5]{Do}.

First of all, recall \cite[Theorem 4.4.1]{Ho2} that any band $B$ decomposes into a semilattice of rectangular bands,
which are the $\D$-classes of $B$. Thus a $\D$-class $D$ of $S$ can be viewed as an $I\times J$ `table' of idempotents
$e_{ij}$ ($i\in I$, $j\in J$), where $\{ R_i\::\: i\in I\}$ and $\{L_j\::\: j\in J\}$ are the $\R$- and $\L$-classes in
$D$ respectively. For $i,k\in I$ and $j,l\in J$ we refer to the tuple $(e_{ij}, e_{il},e_{kj},e_{kl})$ as the
$(i,k;j,l)$ \emph{square}.

Suppose now we have an element $f\in B$ belonging to a $\D$-class above $D$. From the basic theory of bands (see, for
example, \cite[Section 4.4]{Ho2}) we know that $f$ induces idempotent mappings $\sigma\::\: I\rightarrow I$, $i\mapsto
\sigma(i)$, and $\tau\::\: J\rightarrow J$, $j\mapsto (j)\tau$, such that for all $i\in I$, $j\in J$ we have
\[
fe_{ij}=e_{\sigma(i),j},\ e_{ij}f=e_{i,(j)\tau}.
\]
We say that the square $(i,k;j,l)$ is \emph{singular} induced by $f$ if one of the following holds:
\begin{itemize}
\item[\rm (a)] $\sigma(i)=i$, $\sigma(k)=k$ and $(j)\tau=(l)\tau\in\{j,l\}$; or
\item[\rm (b)] $\sigma(i)=\sigma(k)\in\{i,k\}$ and $(j)\tau=j$, $(l)\tau=l$.
\end{itemize}
We talk of a \emph{left-right} or \emph{up-down} singular square depending on whether (a) or (b) applies.

With the above conventions the general presentation we need is as follows:

\begin{pro}[\cite{GR1,Do}]\label{pro2}
The maximal subgroup $H$ of $\ig{B}$ containing $e_{11}\in D$ is presented by
\begin{align}
\label{rel1}
\langle f_{ij}\ (i\in I,\ j\in J) &| && f_{i1}=f_{1j}=1 && (i\in I,\ j\in J),\\
\label{rel2} &&& f_{ij}^{-1}f_{il}=f_{kj}^{-1}f_{kl} && ( (i,k;j,l)\ \text{\rm a singular square in $D$})\rangle.
\end{align}
\end{pro}

\section{Construction of $B_G$}

Let $G$ be any group. Let us choose and fix a presentation $\pre{A}{R}$ for $G$ in which every relation has the form
$ab=c$ for some $a,b,c\in A$. It is clear that $G$ has such a presentation -- for instance the Cayley table would do.
What is less obvious, but nonetheless still true, is that if $G$ is finitely presented then it has a \emph{finite}
presentation of this form. One way of seeing this is as follows: A relation $a_1\dots a_k=b_1\dots b_l$ can be replaced
by two relations of the form $a_1\dots a_k=c$, $b_1\dots b_l=c$, at the expense of introducing a new generator $c$.
Furthermore, the relation $a_1\dots a_k=c$ can be replaced by $k-1$ relations $a_1a_2=d_2$, $d_2a_3=d_3$, \dots,
$d_{k-2}a_{k-1}=d_{k-1}$, $d_{k-1}a_k=c$ of the desired form, with new generators $d_2,\dots,d_{k-1}$.

Define sets
\[
A_0=A\cup\{0\},\ A_0^\prime=\{a^\prime\::\: a\in A_0\},\ I=A_0\cup A_0^\prime,\ J=A_0\cup\{\infty\},
\]
where $0$, $\infty$ and $a^\prime$ ($a\in A_0$) are symbols distinct from each other and those already in $A$. Consider
the semigroup $\T=\T_I^{(l)}\times \T_J^{(r)},$ where $\T_I^{(l)}$ (respectively $\T_J^{(r)}$) is the semigroup of all
mappings $I\rightarrow I$ (resp.\ $J\rightarrow J$) written on the left (resp. right). The semigroup $\T$ has a unique
minimal ideal $K$ consisting of all $(\sigma,\tau)$ with both $\sigma$ and $\tau$ constant. This ideal is naturally
isomorphic to the rectangular band $I\times J$, and we will identify the two. We will visualise $K$ as in Figure 1.

\begin{figure}[ht]
\begin{center}
\begin{pspicture}(30,45)
\psframe[linewidth=0.6](5,0)(30,40) \qline(5,35)(30,35) \psline[linewidth=0.4]{-}(5,20)(30,20) \qline(5,15)(30,15)
\qline(10,0)(10,40) \qline(25,0)(25,40) \put(1,36){$0$} \put(1,16){$0^\prime$} \put(6.5,41){$0$} \put(26,41){$\infty$}
\end{pspicture}
\end{center}
\caption{A visual representation of $K=I\times J$, highlighting the partition $I=A_0\cup A_0^\prime$, as well as the
four distinguished rows and columns.} \label{fig1}
\end{figure}
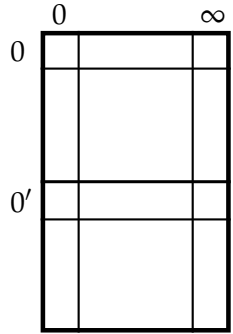

We now define a set $L\subseteq \T$. All elements $(\sigma,\tau)\in L$ will have
\begin{equation}
\label{eq1} \sigma^2=\sigma,\ \tau^2=\tau,\ \ker(\sigma)=\{A_0,A_0^\prime\},\ \im(\tau)=A_0.
\end{equation}
Recall that $\ker(\sigma)$  is the equivalence on $I$ defined by $(i,i')\in\ker(\sigma)$ if and only if
$\sigma(i)=\sigma(i')$, and that it can be identified with the resulting partition of $I$ into equivalence classes.
Therefore, each $(\sigma,\tau)$ will be uniquely determined by $\im(\sigma)$ which must be a two-element set that is a
cross-section of $\{A_0,A_0^\prime\}$, and the value $(\infty)\tau\in A_0$. The elements of $L$ come in four groups:
$\typeZ$ -- the initial pair; $\typeG$, $\typeGbar$ -- the elements arising from the generators $A$; $\typeR$ -- the
elements arising from the relations $R$:

\medskip

\begin{center}
\begin{tabular}{|c|c|c|c|c|}
\hline
 \textbf{Type} & \textbf{Notation} &  \textbf{Indexing} & {\boldmath $\im(\sigma)$} &
{\boldmath $(\infty)\tau$}
\\
\hline\hline $\typeZ$ & $(\sigma_0,\tau_0)$ & -- & $\{0,0^\prime\}$ & $0$
\\
\hline $\typeG$ & $(\sigma_a,\tau_a)$ & $a\in A$ & $\{0,a^\prime\}$ & $a$
\\
\hline $\typeGbar$ & $(\overline{\sigma}_a,\overline{\tau}_a)$ & $a\in A$ & $\{a,a^\prime\}$ & $0$
\\
\hline $\typeR$ & $(\sigma_{\rrel},\tau_{\rrel})$ & $\rrel=(ab,c)\in R $ & $\{b,c^\prime\}$ & $a$
\\
\hline
\end{tabular}
\end{center}

\medskip

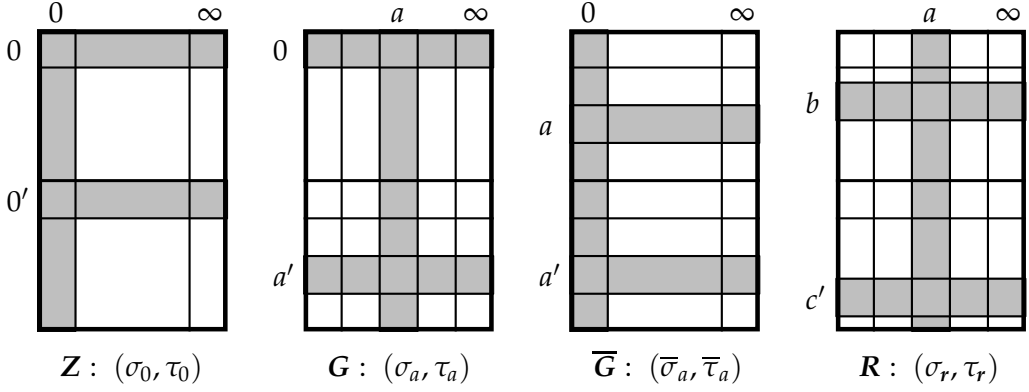
\begin{figure}[ht]
\begin{center}
\begin{pspicture}(0,-7)(140,45)
%
%
\psframe[linewidth=0,fillstyle=solid,fillcolor=lightgray](5,35)(30,40)
\psframe[linewidth=0,fillstyle=solid,fillcolor=lightgray](5,15)(30,20)
\psframe[linewidth=0,fillstyle=solid,fillcolor=lightgray](5,0)(10,40)
\psframe[linewidth=0.6](5,0)(30,40) \qline(5,35)(30,35) \psline[linewidth=0.4]{-}(5,20)(30,20) \qline(5,15)(30,15)
\qline(10,0)(10,40) \qline(25,0)(25,40) \put(1,36){$0$} \put(1,16){$0^\prime$} \put(6.5,41){$0$} \put(26,41){$\infty$}
%
%
\psframe[linewidth=0,fillstyle=solid,fillcolor=lightgray](40,35)(65,40)
\psframe[linewidth=0,fillstyle=solid,fillcolor=lightgray](40,5)(65,10)
\psframe[linewidth=0,fillstyle=solid,fillcolor=lightgray](50,0)(55,40)
\psframe[linewidth=0.6](40,0)(65,40) \qline(40,35)(65,35) \psline[linewidth=0.4]{-}(40,20)(65,20) \qline(40,15)(65,15)
\qline(45,0)(45,40) \qline(60,0)(60,40) \qline(40,5)(65,5) \qline(40,10)(65,10) \qline(50,0)(50,40) \qline(55,0)(55,40)
\put(36,36){$0$} \put(36,6){$a^\prime$} \put(51.5,41){$a$} \put(61,41){$\infty$}
%
%
\psframe[linewidth=0,fillstyle=solid,fillcolor=lightgray](75,25)(100,30)
\psframe[linewidth=0,fillstyle=solid,fillcolor=lightgray](75,5)(100,10)
\psframe[linewidth=0,fillstyle=solid,fillcolor=lightgray](75,0)(80,40)
\psframe[linewidth=0.6](75,0)(100,40) \qline(75,35)(100,35) \psline[linewidth=0.4]{-}(75,20)(100,20)
\qline(75,15)(100,15) \qline(80,0)(80,40) \qline(95,0)(95,40) \qline(75,25)(100,25) \qline(75,30)(100,30)
\qline(75,5)(100,5) \qline(75,10)(100,10) \put(71,26){$a$} \put(71,6){$a^\prime$} \put(76.5,41){$0$}
\put(96,41){$\infty$}
%
%
\psframe[linewidth=0,fillstyle=solid,fillcolor=lightgray](110,28)(135,33)
\psframe[linewidth=0,fillstyle=solid,fillcolor=lightgray](110,2)(135,7)
\psframe[linewidth=0,fillstyle=solid,fillcolor=lightgray](120,0)(125,40)
\psframe[linewidth=0.6](110,0)(135,40) \qline(110,35)(135,35) \psline[linewidth=0.4]{-}(110,20)(135,20)
\qline(110,15)(135,15) \qline(115,0)(115,40) \qline(130,0)(130,40) \qline(110,28)(135,28) \qline(110,33)(135,33)
\qline(110,2)(135,2) \qline(110,7)(135,7) \qline(120,0)(120,40) \qline(125,0)(125,40) \put(106,29){$b$}
\put(106,3){$c^\prime$} \put(121.5,41){$a$} \put(131,41){$\infty$}
%
\put(8,-6){$\typeZ :\ (\sigma_0,\tau_0)$} \put(43,-6){$\typeG :\ (\sigma_a,\tau_a)$} \put(78,-6){$\typeGbar :\
(\overline{\sigma}_a,\overline{\tau}_a)$} \put(113,-6){$\typeR :\ (\sigma_{\rrel},\tau_{\rrel})$}
\end{pspicture}
\end{center}
\caption{The elements $(\sigma_0,\tau_0)$, $(\sigma_a,\tau_a)$, $(\overline{\sigma}_a,\overline{\tau}_a)$ ($a\in A$),
$(\sigma_{\rrel},\tau_{\rrel})$ ($\rrel=(ab,c)\in R$) of $L$. For each $(\sigma,\tau)$ shaded are the two rows
corresponding to $\im(\sigma)$ and one column corresponding to $(\infty)\tau$. They all have
$\ker(\sigma)=\{A_0,A_0^\prime\}$ and $\im(\tau)=A_0$.} \label{fig2}
\end{figure}

\noindent These elements can be visualised as shown in Figure \ref{fig2}.

Because $\ker(\sigma)$ and $\im(\tau)$ are the same for all $(\sigma,\tau)\in L$ it follows that $L$ is a left zero
semigroup (i.e. $xy=x$ for all $x,y\in L$). Furthermore, since $K$ is an ideal in $\T$ (i.e. $xy,yx\in K$ for all $x\in
K$, $y\in\T$), the set $B_G=K\cup L$ is a subsemigroup of $\T$. We remark that, strictly speaking, $B_G$ depends not
only on $G$, but crucially on the chosen presentation for $G$.

\section{Proof of Theorem \ref{th1}}

We will now use the presentation given in Proposition \ref{pro2} to compute the maximal subgroup $H$ of $\ig{B_G}$
containing the idempotent $e_0=(0,0)\in K$. Relations \eqref{rel1} in our context read
\begin{equation}
\label{eq6} f_{0j}=f_{i0}=1\ (i\in I, \ j\in J).
\end{equation}
The remaining relations \eqref{rel2} arise from the singular squares induced by the elements of $L$ acting on $K$. Each
up-down singular square is of one of the following forms:
\[
(a_1,a_2;c_1,c_2),\ (a_1^\prime,a_2^\prime,c_1,c_2)\ (a_1,a_2,c_1,c_2\in A_0).
\]
The square $(a_1,a_2;c_1,c_2)$ yields the relation
\begin{equation}
\label{eq7} f_{a_1,c_1}^{-1} f_{a_1,c_2}=f_{a_2,c_1}^{-1} f_{a_2,c_2}\ (a_1,a_2,c_1,c_2\in A_0).
\end{equation}
Putting $a_1=c_1=0$, $a_2=a$, $c_2=c$ and using \eqref{eq6} yields
\begin{equation}
\label{eq8} f_{a,c}=1\ (a,c\in A_0);
\end{equation}
clearly, all the remaining relations \eqref{eq7} are consequences of \eqref{eq8}. Similarly, the squares
$(a_1^\prime,a_2^\prime,c_1,c_2)$ yield the relations
\begin{equation}
\label{eq9} f_{a^\prime,c}=f_{0^\prime,c}\ (a,c\in A_0).
\end{equation}
(Note that we do not necessarily have $f_{0^\prime,c}=1$, and so cannot deduce $f_{a^\prime,c}=1$.)

Turning to the left-right singular squares, each $(\sigma,\tau)\in L$ induces precisely one. Below we list respectively
the squares introduced by $(\sigma_0,\tau_0)$ of type $\typeZ$, $(\sigma_a,\tau_a)$ of type $\typeG$,
$(\overline{\sigma}_a,\overline{\tau}_a)$ of type $\typeGbar$, and $(\sigma_{\rrel},\tau_{\rrel})$ of type $\typeR$,
together with the relations they yield:
\begin{align}
\label{eq10} &(0,0^\prime;0,\infty): && f_{0,0}^{-1} f_{0,\infty}=f_{0^\prime,0}^{-1} f_{0^\prime,\infty} &&
\\
\label{eq11} &(0,a^\prime;a,\infty): && f_{0,a}^{-1} f_{0,\infty} = f_{a^\prime,a}^{-1} f_{a^\prime,\infty} && (a\in A)
\\
\label{eq12} &(a,a^\prime;0,\infty): && f_{a,0}^{-1} f_{a,\infty} = f_{a^\prime,0}^{-1} f_{a^\prime,\infty} && (a\in A)
\\
\label{eq13} &(b,c^\prime;a,\infty): && f_{b,a}^{-1} f_{b,\infty} = f_{c^\prime,a}^{-1} f_{c^\prime,\infty} &&
(\rrel=(ab,c)\in R).
\end{align}
Using the relations \eqref{eq6}, \eqref{eq8}, \eqref{eq9}, we can transform \eqref{eq10}--\eqref{eq13} into:
\begin{align}
\label{eq14} & f_{0^\prime,\infty}=1 &&
\\
\label{eq15} &f_{a^\prime,\infty} = f_{0^\prime,a} && (a\in A)
\\
\label{eq16} &f_{a,\infty}=f_{a^\prime,\infty}=f_{0^\prime,a} && (a\in A)
\\
\label{eq17} &f_{0^\prime,b}=f_{0^\prime,a}^{-1}f_{0^\prime, c} && (\rrel=(ab,c)\in R).
\end{align}

So, the group $H$ is defined by the generators $f_{i,j}$ ($i\in I$, $j\in J$) and relations \eqref{eq6}, \eqref{eq8},
\eqref{eq9}, \eqref{eq14}--\eqref{eq17}. The relations \eqref{eq6}, \eqref{eq8}, \eqref{eq9},
\eqref{eq14}--\eqref{eq16} can be used simply to eliminate all the generators except $f_{0^\prime,a}$ ($a\in A$).
Replacing each symbol $f_{0^\prime,a}$ by the symbol $a$, the remaining relations \eqref{eq17} become
\[
ab=c \ (\rrel=(ab,c)\in R).
\]
In other words, we obtain the original presentation for $G$. This proves that $H\cong G$.

Finally note that if $\pre{A}{R}$ is a finite presentation, the semigroup $B_G$ is also finite, with
\[
|B_G|=(2|A|+2)(|A|+2)+1+2|A|+|R|,
\]
and this completes the proof of our theorem.

\section{An example, two remarks and an open problem}

It may be instructive to follow in a specific example the sequence of Tietze transformations constituting the brunt of
the above proof. Let us take $G=Q_8$, the quaternion group, with the well known Fibonacci $F(2,3)$ presentation (see
\cite[Section 7.3]{Jo}):
\[
\pre{a,b,c}{ab=c,\ bc=a,\ ca=b}.
\]
The dimension of $K$ in this case is $8\times 5$, and Proposition \ref{pro2} gives a presentation in terms of $40$
generators. This is then simplified by a sequence of generator eliminations, using relations \eqref{rel1}, up-down
singular squares, and left-right singular squares induced by the elements of $L$ of types $\typeZ$, $\typeG$,
$\typeGbar$. In the final step further singular squares are revealed, giving back the original presentation.

\begin{figure}[t]
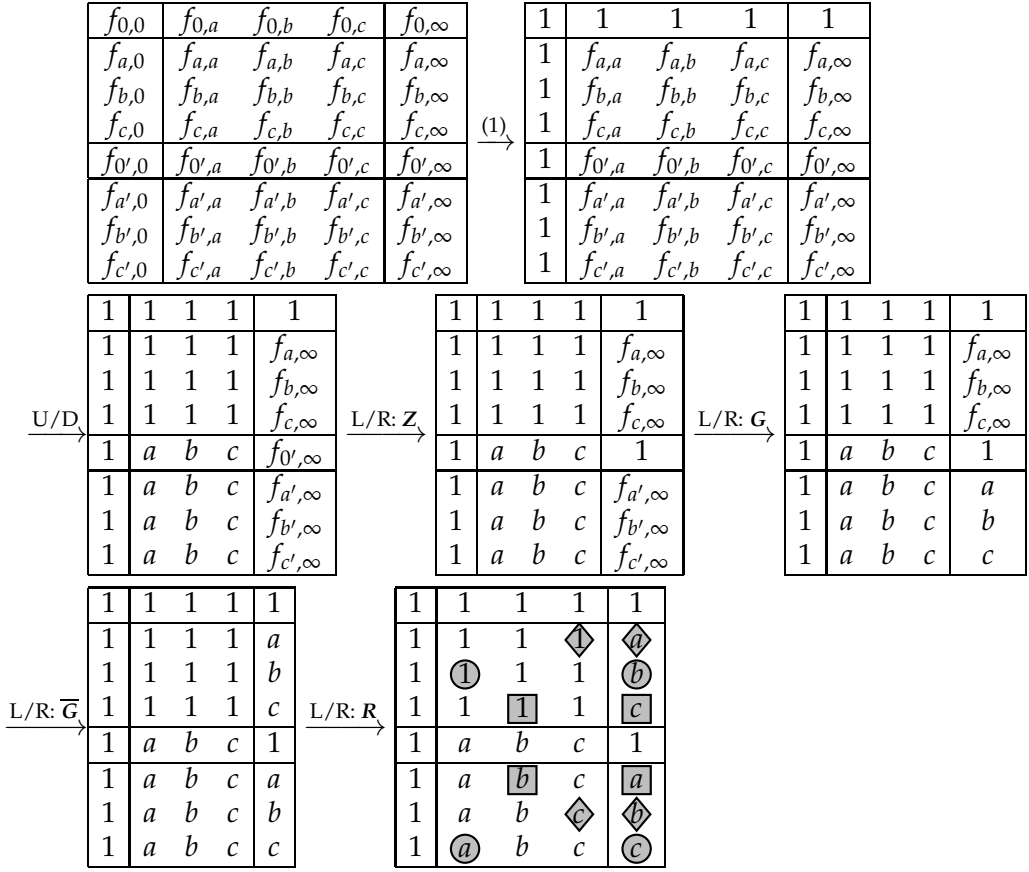

\begin{align*}
&
\begin{array}{|c|ccc|c|}
\hline
f_{0,0} & f_{0,a} & f_{0,b} & f_{0,c} & f_{0,\infty}\\
\hline
f_{a,0} & f_{a,a} & f_{a,b} & f_{a,c} & f_{a,\infty}\\
f_{b,0} & f_{b,a} & f_{b,b} & f_{b,c} & f_{b,\infty}\\
f_{c,0} & f_{c,a} & f_{c,b} & f_{c,c} & f_{c,\infty}\\
\hline
f_{0^\prime,0} & f_{0^\prime,a} & f_{0^\prime,b} & f_{0^\prime,c} & f_{0^\prime,\infty}\\
\hline
f_{a^\prime,0} & f_{a^\prime,a} & f_{a^\prime,b} & f_{a^\prime,c} & f_{a^\prime,\infty}\\
f_{b^\prime,0} & f_{b^\prime,a} & f_{b^\prime,b} & f_{b^\prime,c} & f_{b^\prime,\infty}\\
f_{c^\prime,0} & f_{c^\prime,a} & f_{c^\prime,b} & f_{c^\prime,c} & f_{c^\prime,\infty}\\
\hline
\end{array}
\xrightarrow{\eqref{rel1}}
\begin{array}{|c|ccc|c|}
\hline
1&1&1&1&1\\
\hline
1& f_{a,a} & f_{a,b} & f_{a,c} & f_{a,\infty}\\
1 & f_{b,a} & f_{b,b} & f_{b,c} & f_{b,\infty}\\
1& f_{c,a} & f_{c,b} & f_{c,c} & f_{c,\infty}\\
\hline
1 & f_{0^\prime,a} & f_{0^\prime,b} & f_{0^\prime,c} & f_{0^\prime,\infty}\\
\hline
1 & f_{a^\prime,a} & f_{a^\prime,b} & f_{a^\prime,c} & f_{a^\prime,\infty}\\
1 & f_{b^\prime,a} & f_{b^\prime,b} & f_{b^\prime,c} & f_{b^\prime,\infty}\\
1 & f_{c^\prime,a} & f_{c^\prime,b} & f_{c^\prime,c} & f_{c^\prime,\infty}\\
\hline
\end{array}
\\
\xrightarrow{\text{U/D}} &
\begin{array}{|c|ccc|c|}
\hline
1&1&1&1&1\\
\hline
1& 1&1&1 & f_{a,\infty}\\
1 & 1&1&1 & f_{b,\infty}\\
1& 1&1&1 & f_{c,\infty}\\
\hline
1 & a&b&c & f_{0^\prime,\infty}\\
\hline
1 & a&b&c& f_{a^\prime,\infty}\\
1 & a&b&c& f_{b^\prime,\infty}\\
1 & a&b&c & f_{c^\prime,\infty}\\
\hline
\end{array}
\xrightarrow{\text{L/R: } \typeZ}
\begin{array}{|c|ccc|c|}
\hline
1&1&1&1&1\\
\hline
1& 1&1&1 & f_{a,\infty}\\
1 & 1&1&1 & f_{b,\infty}\\
1& 1&1&1 & f_{c,\infty}\\
\hline
1 & a&b&c& 1\\
\hline
1 & a&b&c& f_{a^\prime,\infty}\\
1 & a&b&c & f_{b^\prime,\infty}\\
1 & a&b&c & f_{c^\prime,\infty}\\
\hline
\end{array}
\xrightarrow{\text{L/R: } \typeG}
\begin{array}{|c|ccc|c|}
\hline
1&1&1&1&1\\
\hline
1& 1&1&1 & f_{a,\infty}\\
1 & 1&1&1 & f_{b,\infty}\\
1& 1&1&1 & f_{c,\infty}\\
\hline
1 & a&b&c& 1\\
\hline
1 & a&b&c& a\\
1 & a&b&c & b\\
1 & a&b&c & c\\
\hline
\end{array}
\\
\xrightarrow{\text{L/R: } \typeGbar} &
\begin{array}{|c|ccc|c|}
\hline
1&1&1&1&1\\
\hline
1& 1&1&1 & a\\
1 & 1&1&1 & b\\
1& 1&1&1 & c\\
\hline
1 & a&b&c& 1\\
\hline
1 & a&b&c& a\\
1 & a&b&c & b\\
1 & a&b&c & c\\
\hline
\end{array}
\xrightarrow{\text{L/R: } \typeR}
\begin{array}{|c|ccc|c|}
\hline
1&1&1&1&1\\
\hline
1& 1&1&\indiamond{$1$} & \indiamond{$a$}\\
1 & \incircle{$1$}&1&1 & \incircle{$b$}\\
1& 1&\insquare{$1$}&1 & \insquare{$c$}\\
\hline
1 & a&b&c& 1\\
\hline
1 & a&\insquare{$b$}&c& \insquare{$a$}\\
1 & a&b&\indiamond{$c$} & \indiamond{$b$}\\
1 & \incircle{$a$}&b&c & \incircle{$c$}\\
\hline
\end{array}
\end{align*}
\caption{The sequence of Tietze transformations constituting the proof of Theorem \ref{th1}.} \label{fig3}

\end{figure}

If we record the original generators in a natural $8\times 5$ grid, this process may be encapsulated as shown in Figure
\ref{fig3}.

\begin{rem}
\label{rem1} It is possible to describe completely the structure of the free idempotent generated semigroup $\ig{B_G}$.
By known results (see e.g.\ \cite[(IG1)--(IG4)]{GR1}) $\ig{B_G}$ has precisely two regular $\D$-classes. The `upper'
one is a left zero semigroup $\overline{L}$ isomorphic to $L$ (as all products in $L$ are basic), while the `lower' one
$\overline{K}$, the completely simple minimal ideal, has a Rees matrix representation with structure group $G$ and
(normalised) sandwich matrix $(a_{ji}^{-1})$, where $(a_{ij})$ is the $|I|\times|J|$ table that is the end-product of
Tietze transformations performed in the proof of Theorem \ref{th1} (in our example this is the last table in Figure
\ref{fig3}). We claim that in fact $\ig{B_G}=\overline{L}\cup\overline{K}$. To confirm this, and see that the structure
is completely determined, we need to show how to write products $ef$ and $fe$ with $e\in L$, $f\in K$ as products of
idempotents from $K$ in $\ig{B}$. For the product $ef$ note that there exists $g\in K$ such that $f\R g$ and $ge=g$;
both pairs $\{e,g\} $ and $\{f,g\}$ are critical and we have $ef=egf=hf$, where $h=eg\in K$. The product $fe$ can be
treated similarly.
\end{rem}

\begin{rem}
\label{rem2} Associated to the biorder $E$ of idempotents of a \emph{regular} semigroup $S$ there is another free
idempotent generated object $\rig{E}$, the \emph{free regular idempotent generated semigroup} on $E$. It is the largest
\emph{regular} semigroup with the biorder of idempotents $E$, and its presentation can be obtained by adding further
relations to the defining presentation for $\ig{E}$. For definition and references we refer the reader to \cite{GR1}.
In particular, from (IG1)--(IG4), (RIG1), (RIG2) in \cite{GR1} it follows that $\ig{B_G}=\rig{B_G}$.
\end{rem}

One way of interpreting Remarks \ref{rem1}, \ref{rem2} is to say that the word problem for $\ig{B_G}$ is decidable if
and only if the word problem for $G$ is decidable. It is the authors' belief that the next stage in the ongoing
exploration of free idempotent generated semigroups is precisely an analysis of the word problem for $\ig{S}$. This at
present seems a daunting task, even in the case where $S$ is finite. Nonetheless, we propose the following problem
which may just be within reach at this stage:


\begin{que}
Let $B$ be a finite band such that all maximal subgroups of $\ig{B}$ have recursively soluble word problems. Is the
word problem of $\ig{B}$ necessarily recursively soluble?
\end{que}

\begin{ack}
The research of the first author is supported by the Ministry of Education, Science and Technological Development of
the Republic of Serbia through Grant No.174019, and by a grant (Contract 114--451--2675/2012) of the Secretariat of
Science and Technological Development of the Autonomous Province of Vojvodina. Also, the first author gratefully
acknowledges the hospitality of the School of Mathematics and Statistics of the University of St Andrews, where this
research was carried out.
\end{ack}



\begin{thebibliography}{99}
\frenchspacing

\bibitem{BMM1}
M. Brittenham, S. W. Margolis and J. Meakin, Subgroups of free idempotent generated semigroups need not be free,
\emph{J. Algebra} \textbf{321} (2009), 3026--3042.

\bibitem{Do}
I. Dolinka, A note on maximal subgroups of free idempotent generated semigroups over bands, \emph{Periodica Math.
Hungar.} \textbf{65} (2012), 97--105.

\bibitem{DG}
I. Dolinka and R. Gray, Maximal subgroups of free idempotent generated semigroups over the full linear monoid,
\emph{Trans. Amer. Math. Soc.}, to appear. \url{arXiv: 1112.0893}

\bibitem{E1}
D. Easdown, Biordered sets of bands, \emph{Semigroup Forum} \textbf{29} (1984), 241--246.

\bibitem{E2}
D. Easdown, Biordered sets are biordered subsets of idempotents of semigroups, \emph{J. Austral. Math. Soc. Ser. A}
\textbf{37} (1984), 258--268.

\bibitem{E3}
D. Easdown, Biordered sets come from semigroups, \emph{J. Algebra} \textbf{96} (1985), 581--591.

\bibitem{Er}
J. A. Erdos, On products of idempotent matrices, \emph{Glasgow Math. J.} \textbf{8} (1967), 118--122.

\bibitem{GY}
V. Gould, D. Yang,  Every group is the maximal subgroup of a naturally occurring free idempotent generated semigroup,
\url{arXiv:1209.1242}

\bibitem{GR1}
R. Gray and N. Ru\v skuc, On maximal subgroups of free idempotent generated semigroups, \emph{Israel J. Math.}
\textbf{189} (2012), 147--176.

\bibitem{GR2}
R. Gray and N. Ru\v skuc, Maximal subgroups of free idempotent generated semigroups over the full transformation
monoid, \emph{Proc. London Math. Soc.} \textbf{105} (2012), 997--1018.

\bibitem{Hi}
P. M. Higgins, \emph{Techniques of Semigroup Theory}, Oxford University Press, New York, 1992.

\bibitem{Ho1}
J. M. Howie, The subsemigroup generated by the idempotents of a full transformation semigroup, \emph{J. London Math.
Soc.} \textbf{41} (1966), 707--716.

\bibitem{Ho2}
J. M. Howie, \emph{Fundamentals of Semigroup Theory}, Oxford University Press, New York, 1995.

\bibitem{Jo}
D. L. Johnson, \emph{Presentations of Groups}, LMS Student Texts Vol. 15, Cambridge University Press, Cambridge, 1990.

\bibitem{McE}
B. McElwee, Subgroups of the free semigroup on a biordered set in which principal ideals are singletons, \emph{Comm.
Algebra} \textbf{30} (2002), 5513--5519.

\bibitem{Na}
K. S. S. Nambooripad, Structure of regular semigroups. I, \emph{Mem. Amer. Math. Soc.} \textbf{22} (1979), no. 224,
vii+119 pp.

\bibitem{NP}
K. S. S. Nambooripad and F. Pastijn, Subgroups of free idempotent generated regular semigroups, \emph{Semigroup Forum}
\textbf{21} (1980), 1--7.

\bibitem{Ru}
N. Ru\v skuc, Presentations for subgroups of monoids, \emph{J. Algebra} \textbf{220} (1999), 365--380.

\bibitem{Pa}
F. Pastijn, The biorder on the partial groupoid of idempotents of a semigroup, \emph{J. Algebra} \textbf{65} (1980),
147--187.

\bibitem{Pe}
M. Petrich, \emph{Introduction to Semigroups}, Merrill, Columbus, 1973.

\end{thebibliography}
\end{document}